\documentclass[11pt]{amsart}
\usepackage{amscd}      
\usepackage{eucal,amsfonts,graphics,epic,amsbsy, stmaryrd}
\usepackage{amssymb}
\usepackage{amsmath}
\usepackage{xypic}      
\LaTeXdiagrams          
\usepackage[all]{xy}

\usepackage{latexsym}
\usepackage{epsfig}
\usepackage{amsfonts}
\usepackage{enumerate}

\newtheorem{prop}{Proposition}[section]

\newtheorem{thm}[prop]{Theorem}

\theoremstyle{definition}

\newtheorem{rmk}[prop]{Remark}

\newtheorem{example}[prop]{Example}

            {\nolinebreak $\Box$ \end{trivlist}}

\newcommand{\noprint}[1]{}

\newcommand{\XX}{{\mathfrak X}}

\newcommand{\YY}{{\mathfrak Y}}

\newcommand{\cc}{{\mathbb C}}

\newcommand{\ldiag}[1]%
       {\makebox[0cm]{${\scriptstyle#1}\downarrow\phantom{\scriptstyle#1}$}}
\newcommand{\ldiagup}[1]%
       {\makebox[0cm]{${\scriptstyle#1}\uparrow\phantom{\scriptstyle#1}$}}
\newcommand{\rdiag}[1]%
       {\makebox[0cm]{$\phantom{\scriptstyle#1}\downarrow{\scriptstyle#1}$}}
\newcommand{\sediagr}[1]%
       {\makebox[0cm]{$\phantom{\scriptstyle#1}\searrow{\scriptstyle#1}$}}
\newcommand{\nediagr}[1]%
       {\makebox[0cm]{$\phantom{\scriptstyle#1}\nearrow{\scriptstyle#1}$}}
\newcommand{\rdiagup}[1]%
       {\makebox[0cm]{$\phantom{\scriptstyle#1}\uparrow{\scriptstyle#1}$}}
\newcommand{\swdiag}[1]%
       {\makebox[0cm]{$\phantom{\scriptstyle#1}\swarrow{\scriptstyle#1}$}}
\newcommand{\sediag}[1]%
       {\makebox[0cm]{${\scriptstyle#1}\searrow\phantom{\scriptstyle#1}$}}
\newcommand{\nediag}[1]%
       {\makebox[0cm]{${\scriptstyle#1}\nearrow\phantom{\scriptstyle#1}$}}

\newcommand{\doublearrowstack}[2]%
                      {{{{\scriptstyle#1}\atop{\textstyle\longrightarrow}}\atop{{\textstyle\longrightarrow}\atop{\scriptstyle#2}}}}
\newcommand{\rightleftarrowstack}[2]%
                      {{{{\scriptstyle#1}\atop{\textstyle\longrightarrow}}\atop{{\textstyle\longleftarrow}\atop{\scriptstyle#2}}}}
\newcommand{\leftrightarrowstack}[2]%
                      {{{{\scriptstyle#1}\atop{\textstyle\longleftarrow}}\atop{{\textstyle\longrightarrow}\atop{\scriptstyle#2}}}}

\newcommand{\overtoparrow}%
{\makebox[0cm]{\beginpicture \setcoordinatesystem units
<.8cm,.4cm> point at 0 0 \setplotarea x from -3 to 3, y from 0 to
1 \setquadratic \plot -3 0 0 1 3 0 / \put{\vector(3,-1){0}}[Bl] at
3 0
\endpicture}}

\newcommand{\underbottomarrow}%
{\makebox[0cm]{\beginpicture \setcoordinatesystem units
<.8cm,.4cm> point at 0 0 \setplotarea x from -3 to 3, y from 0 to
1 \setquadratic \plot -3 1 0 0 3 1 / \put{\vector(3,1){0}}[Bl] at
3 1
\endpicture}}

\newcommand{\ses}[5]%
{0\longrightarrow#1\stackrel{#2}{ \longrightarrow}#3\stackrel{#4}{
\longrightarrow}#5\longrightarrow0}

\newcommand{\dt}[6]%
{#1\stackrel{#2}{longrightarrow}#3
\stackrel{#4}{\longrightarrow}#5 \stackrel{#6}{\longrightarrow}
#1[1]}

\newcommand{\cat}[1]%
{(\mbox{\rm #1})}

 \def\bbZ{{\mathbb Z}}

\newcommand{\clG}{{\mathcal{G}}}

\newcommand{\clK}{{\mathcal{K}}}

\title{On Gromov-Witten theory of root gerbes}
\date{\today}
\author{Elena Andreini}
\address{Max-Planck-Institut f\"ur Mathematik\\Vivatsgasse 7\\
53111 Bonn\\ Germany}
\email{andreini.elena@gmail.com}
\author{Yunfeng Jiang}
\address{Department of Mathematics\\ University of Utah\\ 155 S 1400 E JWB 233\\Salt Lake city\\ UT 84112\\ USA}
\email{jiangyf@math.utah.edu}
\author{Hsian-Hua Tseng}
\address{Department of Mathematics\\ University of Wisconsin-Madison\\ Van Vleck Hall, 480 Lincoln Drive \\Madison\\ WI 53706-1388 \\ USA}
\email{tseng@math.wisc.edu}

\begin{document}
\sloppy \maketitle
\begin{abstract}
This research announcement discusses our results on Gromov-Witten theory of root gerbes. 
A complete calculation  of genus $0$ Gromov-Witten theory of $\mu_{r}$-root gerbes over a smooth base scheme is obtained by a direct analysis of virtual fundamental classes. Our result verifies the genus $0$ part of the so-called decomposition conjecture which compares Gromov-Witten theory of \'etale gerbes with that of the bases. We also verify this conjecture in all genera for toric gerbes over toric Deligne-Mumford stacks. 
\end{abstract}


\maketitle

\section{Introduction}

Orbifold Gromov-Witten theory, constructed in symplectic category by Chen-Ruan \cite{CR2}  and in algebraic category by Abramovich, Graber and Vistoli \cite{AGV1}, \cite{AGV2}, has been an area of active research in recent years. Calculations of orbifold Gromov-Witten invariants in examples present numerous new challenges, see \cite{CCLT}, \cite{CCIT}, \cite{ts}, and \cite{BC} for examples. 

\'Etale gerbes over a smooth base are examples of smooth Deligne-Mumford stacks. Let $\XX$ be a smooth Deligne-Mumford stack and  $G$  a finite group  scheme over $\XX$. Intuitively one can think of a $G$-banded gerbe over  $\XX$ as a fibre bundle over $\XX$  with fibre the classifying stack $\mathcal{B}G$. A detailed definition of gerbes can be found in, for example, \cite{Gir}, \cite{Bree94}, \cite{ehkv}. Our main goal is to compute Gromov-Witten theory of $G$-banded gerbes. 

For trivial $G$-banded gerbes $\XX\times \mathcal{B}G$, the computation of their Gromov-Witten invariants is handled as a special case of a general product formula for orbifold Gromov-Witten invariants of product stacks $\XX\times\YY$. We will not discuss this here, see \cite{AJT} for more details.
  
Root gerbes provide an interesting class of non-trivial gerbes. Let $\XX$ be a smooth Deligne-Mumford stack.  Let $\mathcal{L}\rightarrow \XX$ be a line bundle over $\XX$ and $r$ a positive integer. The stack $\sqrt[r]{\mathcal{L}/\XX}$ of $r$-th roots of $\mathcal{L}$ is  a smooth Deligne-Mumford stack, and is a $\mu_{r}$-gerbe over $\XX$. Our study of Gromov-Witten theory of root gerbes is aided by the so-called {\em decomposition conjecture} \cite{HHPS} in physics, which states\footnote{One can also formulate the decomposition conjecture for arbitrary $G$-gerbes (which is more general than a $G$-banded gerbe). The conjecture states that Gromov-Witten theory of the $G$-gerbe is equivalent to certain twist of the Gromov-Witten theory of some \'etale cover of the base. Details of the conjecture in this generality will be discussed elsewhere.} in mathematical terms that Gromov-Witten theory of a root gerbe $\sqrt[r]{\mathcal{L}/\XX}$ over $\XX$ should be equivalent to Gromov-Witten theory of the disjoint union of $r$ copies of $\XX$ after a change of variables. 

In this note, we present results on computations of genus $0$ Gromov-Witten theory of root gerbes over smooth varieties.  Our results verifies the decomposition conjecture in genus $0$ for these root gerbes. We also discuss the case of {\em toric gerbes}\footnote{Note that these gerbes are iterated root gerbes over  toric Deligne-Mumford stacks \cite{FMN}, \cite{JT}.}, for which we verify the decomposition conjecture in all genera by applying some sophisticated techniques in toric Gromov-Witten theory. Detailed proofs of results discussed in this note will be given in \cite{AJT2}.

\subsection*{Acknowledgments} We thank D. Abramovich, A. Bayer,  K. Behrend,  B. Fantechi,  P. Johnson, A. Kresch,  Y. Ruan and A. Vistoli for valuable discussions. H.-H. T. is grateful to T. Coates, A. Corti, H. Iritani, and X. Tang for related collaborations. H.-H. T. is supported in part by NSF grant DMS-0757722.

\section{Results on root gerbes}

\subsection{Root gerbes}
Let $\XX$ be a proper smooth Deligne-Mumford stack,  and $\mathcal{L}$ a line bundle over $\XX$. Let $r>0$ be an integer. 
Recall \cite{AGV2}, \cite{Ca}  that the stack $\sqrt[r]{\mathcal{L}/\XX}$ of $r$-th roots of $\mathcal{L}$ is defined to be the $\XX$-groupoid  whose objects
over an  $\XX$-scheme  $f: Y\to \XX$ are pairs $(M,\varphi)$, with $M$ a line bundle over
$Y$ and $\varphi: M^{\otimes r}\to f^{*}\mathcal{L}$ an
isomorphism. An arrow from $(M,\varphi)$ to $(N,\psi)$ lying over
a $\XX$-morphism $h: (Y,f)\to (Z,g)$ is an isomorphism
$\rho: M\to h^{*}N$ such that
$\varphi=(h^{*}\psi)\circ \rho^{\otimes r}$:
\[\xymatrix{
M^{\otimes r} \ar[r]^{\rho^{\otimes r}} \ar[d]_{\varphi} & h^{*}N^{\otimes r} \ar[d]^{h^{*}\psi} \\
f^{*}\mathcal{L}  \ar[r]^{\cong}  & h^{*}g^{*}\mathcal{L} .}
\]
Alternatively, the stack $\sqrt[r]{\mathcal{L}/\XX}$ may be presented as the quotient stack
$[\mathcal{L}^{*}/\mathbb{C}^{*}]$, where $\mathcal{L}^{*}$ is the complement of the zero section
inside the total space of $\mathcal{L}$, and $\mathbb{C}^{*}$ acts via $\lambda\cdot z:= \lambda^rz$.
Clearly $\sqrt[r]{\mathcal{L}/\XX}$ is a smooth Deligne-Mumford stack. It is also easy to see that $\sqrt[r]{\mathcal{L}/\XX}\to \XX$ is a $\mu_r$-banded gerbe. See \cite{AGV2} and \cite{Ca} for more discussions.

In this note we are primarily concerned with the case that the base $\XX$ is a  smooth projective variety\footnote{In this case, the stack $\sqrt[r]{\mathcal{L}/X}$ is a toric stack bundle introduced in \cite{Jiang}.}
 $X$. 

By definition the Chen-Ruan orbifold cohomology groups (see \cite{CR1}, \cite{AGV1} for more details) are cohomology groups of the {\em inertia stacks}. Since  the root gerbe $\sqrt[r]{\mathcal{L}/X}$ is naturally banded by $\mu_r$  its  inertia stack $I\sqrt[r]{\mathcal{L}/X}$ is a disjoint union of $r$ $\mu_r$-gerbes over $X$. Thus the Chen-Ruan orbifold cohomology  with rational coefficients $H_{CR}^{*}(\sqrt[r]{\mathcal{L}/X}, \mathbb{Q})$ is isomorphic to a direct sum of $r$ copies of $H^*(X, \mathbb{Q})$.

\subsection{Moduli of stable maps}

One of our results is the construction of moduli stack of twisted stable maps. Let $\beta\in H_{2}(X,\bbZ)$ and let  $\mathcal{K}_{0,n}(\sqrt[r]{\mathcal{L}/X},\beta)$ be the moduli stack of genus zero, $n$-point twisted stable maps to $\sqrt[r]{\mathcal{L}/X}$ of degree $\beta$ in the sense of Abramovich-Graber-Vistoli \cite{AGV2}. Given $[f:\mathcal{C}\to \sqrt[r]{\mathcal{L}/X}]\in \mathcal{K}_{0,n}(\sqrt[r]{\mathcal{L}/X},\beta)$, the stack structures at marked points of $\mathcal{C}$ determine an $n$-tuple $\vec{g}=(g_1,\cdots,g_n)$ of elements in $\mu_{r}$. We may write $g_{i}=e^{2\pi i m_{i}/r_{i}}$, where  $0\leq m_i\leq r_{i}-1$ and $(m_i, r_i)=1$. By Riemann-Roch for twisted curves, $\vec{g}$ satisfies the condition 
$$\prod_{1\leq i\leq n}g_{i}=e^{2\pi i k/r}, \text{ where }k=\int_\beta c_{1}(\mathcal{L}).$$
Such $n$-tuples $\vec{g}$ are called {\em admissible vectors}.

Let $$\mathcal{K}_{0,n}(\sqrt[r]{\mathcal{L}/X},\beta)^{\vec{g}}\subset \mathcal{K}_{0,n}(\sqrt[r]{\mathcal{L}/X},\beta)$$ be the locus which parametrizes stable maps with admissible vector $\vec{g}$. Post-composition with the natural map $\sqrt[r]{\mathcal{L}/X}\to X$ gives a map $$p: \mathcal{K}_{0,n}(\sqrt[r]{\mathcal{L}/X},\beta)^{\vec{g}}\to \overline{\mathcal{M}}_{0,n}(X, \beta),$$
where $\overline{\mathcal{M}}_{0,n}(X, \beta)$ is the moduli stack of genus zero, $n$-point degree $\beta$ stable maps to $X$.
We prove  a structure result for the moduli stack 
$\mathcal{K}_{0,n}(\sqrt[r]{\mathcal{L}/X},\beta)^{\vec{g}}$. More precisely,  we obtain a construction of $\mathcal{K}_{0,n}(\sqrt[r]{\mathcal{L}/X},\beta)^{\vec{g}}$ as follows.
\begin{enumerate}
\item
construct a stack $P^{\vec{g}}_n$ over $\overline{\mathcal{M}}_{0,n}(X, \beta)$ as a category fibered in groupoids whose objects are certain morphisms of  logarithmic  structures \cite{KatoLog}, \cite{MO}. This construction, which can be seen as a generalization of the root construction,  introduces additional  automorphisms along the locus in $\overline{\mathcal{M}}_{0,n}(X, \beta)$ parametrizing maps with singular domains;

\item
show that $\clK_{0,n}(\clG,\beta)^{\vec{g}}$ is a $\mu_r$-banded gerbe  over $P^{\vec{g}}_n$ and therefore  $p$ factors through $P^{\vec{g}}_n$.
\end{enumerate}

The $\mu_r$-gerbe in step (2) is the root gerbe of a line bundle which we can explicitly write down. Our result thus generalizes a result of Bayer and Cadman \cite{BC} on the moduli stack of twisted stable maps to the classifying stack $\mathcal{B}\mu_r$.

Our strategy of computing genus $0$ Gromov-Witten invariants of $\sqrt[r]{\mathcal{L}/X}$ is to relate them with invariants of $X$. To do so, we carry out a comparison with respect to the map $p$  between the perfect relative obstruction theories on $\mathcal{K}_{0,n}(\sqrt[r]{\mathcal{L}/X},\beta)^{\vec{g}}$ and $\overline{\mathcal{M}}_{0,n}(X,\beta)$  and deduce from there the following push-forward formula for virtual fundamental classes:

\begin{thm}\label{pushforward_vir}
$$p_{*}[\mathcal{K}_{0,n}(\sqrt[r]{\mathcal{L}/X},\beta)^{\vec{g}}]^{vir}= 
 \frac{1}{r}[\overline{\mathcal{M}}_{0,n}(X,\beta)]^{vir}.$$
\end{thm} 

\subsection{Gromov-Witten invariants}
Let
$$\pi: \clG=\sqrt[r]{\mathcal{L}/X} \to X$$ 
be a $\mu_{r}$-root gerbe.
Then the inertia stack admits the following decomposition
$$I\clG=\coprod_{g\in \mu_{r}}\clG_{g},$$
where $\clG_{g}$ is a  root gerbe isomorphic to  $\clG$. Let $\pi_{g}:\clG_{g}\to X$
be the induced morphism.  On each component there is an isomorphism between the rational cohomology groups 
$$\pi_{g}^{*}: H^{*}(X,\mathbb{Q})\overset{\simeq}{\longrightarrow} H^{*}(\clG_{g}, \mathbb{Q}).$$

Let $\vec{g}=(g_1,...,g_n)$ be an admissible vector. There are evaluation maps $$ev_i: \mathcal{K}_{0,n}(\sqrt[r]{\mathcal{L}/X},\beta)^{\vec{g}}\to \clG^{rig}_{g_i},$$
where $\clG^{rig}_{g_i}$ is a component of  the  {\em rigidified inertia stack} $I\clG^{rig}=\cup_{g\in \mu_r} \clG^{rig}_g$  (see \cite{AGV2}, Section 4.4 for the  definition). Although the evaluation maps  $ev_i$ do not take values in  $I\clG$, as explained in \cite{AGV2}, Section 6.1.3, one can still define a pull-back map at cohomology level,
$$ev_i^*: H^*(\clG_{g_i}, \mathbb{Q})\to H^*(\mathcal{K}_{0,n}(\sqrt[r]{\mathcal{L}/X},\beta)^{\vec{g}}, \mathbb{Q}).$$
Given $\delta_i\in H^*(\clG_{g_{i}},\mathbb{Q})$ for $1\leq i\leq n$ and   integers $k_i\geq 0,1\leq i\leq n$,  one can define descendant orbifold Gromov-Witten invariants 
$$\langle \delta_1\bar{\psi}_1^{k_1},\cdots, \delta_n\bar{\psi}_n^{k_n}\rangle_{0, n, \beta}^\clG:=\int_{[\mathcal{K}_{0,n}(\sqrt[r]{\mathcal{L}/X},\beta)^{\vec{g}}]^{vir}}\prod_{i=1}^n ev_i^*(\delta_i)\bar{\psi}_i^{k_i},$$
where $\overline{\psi}_i$ are the pullback of the first  Chern classes  of the tautological line bundles  over $\overline{\mathcal{M}}_{g,n}(X,\beta)$. See \cite{AGV2} for more discussion on descendant classes.

For  classes $\delta_i\in H^*(\clG_{g_i},\mathbb{Q})$, 
set $\overline{\delta}_i=(\pi_{g_{i}}^{*})^{-1}(\delta_i)$. Descendant 
Gromov-Witten invariants  $\langle\overline{ \delta}_1\bar{\psi}_1^{k_1},\cdots,\overline{\delta}_n\bar{\psi}_n^{k_n}\rangle_{0, n, \beta}^X$ of $X$ are similarly defined. 
Theorem \ref{pushforward_vir} implies
\begin{thm}\label{GW-inv1}
\begin{equation*}
 \langle \delta_1\bar{\psi}_1^{k_1},..., \delta_n\bar{\psi}_n^{k_n}\rangle_{0,n,\beta}^\clG=\frac{1}{r} \langle \overline{\delta}_1\bar{\psi}_1^{k_1},\cdots,\overline{\delta}_n\bar{\psi}_n^{k_n}\rangle_{0,n,\beta}^X.
\end{equation*}
\end{thm}
Moreover, if $\vec{g}$ is not admissible, then the Gromov-Witten invariants vanish.

In the following we use complex number $\cc$ as coefficients for the cohomology. For $\overline{\alpha}\in H^*(X, \mathbb{C})$ and an irreducible representation $\rho$ of $\mu_r$, we define $$\overline{\alpha}_\rho:=\frac{1}{r}\sum_{g\in \mu_r} \chi_{\rho}(g^{-1})\pi_g^*(\overline{\alpha}),$$ where $\chi_\rho$ is the character of $\rho$. The map $(\overline{\alpha},\rho)\mapsto \overline{\alpha}_\rho$  clearly defines an additive isomorphism
$$\bigoplus_{[\rho]\in \widehat{\mu_r}}H^*(X)_{[\rho]}\simeq H^*(I\clG, \mathbb{C}),$$ where $\widehat{\mu_r}$ is the set of isomorphism classes of irreducible representations of $\mu_r$, and for $[\rho]\in\widehat{\mu_r}$ we define $H^*(X)_{[\rho]}:= H^*(X, \mathbb{C})$.

Theorem \ref{GW-inv1} together with orthogonality relations of characters of  
$\mu_r$ implies the following 

\begin{thm}\label{GW-inv2}
\begin{equation*}
\begin{split}
&\langle \overline{\alpha}_{1\rho_1}\bar{\psi}_1^{k_1},..., \overline{\alpha}_{n\rho_n}\bar{\psi}_n^{k_n} \rangle_{0,n,\beta}^{\clG}\\
=&\begin{cases}
\frac{1}{r^2}\langle \overline{\alpha}_1\bar{\psi}_1^{k_1},\cdots,\overline{\alpha}_n\bar{\psi}_n^{k_n}\rangle_{0,n,\beta}^X\chi_\rho(\exp(\frac{-2\pi i \int_\beta c_1(\mathcal{L})}{r})) &\text{if }\rho_1=\rho_2=...=\rho_n=:\rho,\\
0&\text{otherwise}\,.
\end{cases}
\end{split}
\end{equation*}
\end{thm}

We may reformulate this in terms of generating functions. Let $$\{\overline{\phi}_i\, |\, 1\leq i\leq \text{rank}H^*(X, \mathbb{C})\}\subset H^*(X, \mathbb{C})$$ be an additive basis. According to the discussion above, the set $$\{\overline{\phi}_{i\rho}\,|\, 1\leq i\leq \text{rank}H^*(X, \mathbb{C}), [\rho]\in \widehat{\mu_r}\}$$ is an additive basis of $H^*(I\clG, \mathbb{C})$. Recall that the genus $0$ descendant potential of $\clG$ is defined to be 
$$\mathcal{F}^0_{\clG}(\{t_{i\rho, j}\}_{1\leq i\leq \text{rank}H^{*}(X,\mathbb{C}), \rho\in \widehat{\mu_r}, j\geq 0}; Q):=\sum_{\overset{n\geq 0, \beta\in H_2(X,\mathbb{Z})}{i_1,...,i_n; \rho_1,...,\rho_n; j_1,...,j_n}}\frac{Q^\beta}{n!}\prod_{k=1}^nt_{i_k\rho_k, j_k}\langle\prod_{k=1}^n \overline{\phi}_{i_k\rho_k}\bar{\psi}_k^{j_k} \rangle_{0,n,\beta}^{\clG}.$$
The descendant potential $\mathcal{F}^0_{\clG}$ is a formal power series in variables $t_{i\rho, j}, 1\leq i\leq \text{rank}H^*(X,\mathbb{C}), \rho\in \widehat{\mu_r}, j\geq 0$ with coefficients in the Novikov ring $\mathbb{Q}[[\overline{NE}(X)]]$, where $\overline{NE}(X)$ is the effective Mori cone of the coarse moduli space of $\clG$. Here $Q^\beta$ are formal variables labeled by classes $\beta\in \overline{NE}(X)$. See e.g. \cite{ts} for more discussion on descendant potentials for orbifold Gromov-Witten theory. 

Similarly the genus $0$ descendant potential of $X$ is defined to be 
$$\mathcal{F}^0_X(\{t_{i,j}\}_{1\leq i\leq \text{rank}H^*(X,\mathbb{C}), j\geq 0};Q):=\sum_{\overset{n\geq 0, \beta\in H_2(X, \mathbb{Z})}{i_1,...,i_n; j_1,...,j_n}}\frac{Q^\beta}{n!}\prod_{k=1}^nt_{i_k, j_k}\langle \prod_{k=0}^n \overline{\phi}_{i_k}\bar{\psi}_k^{j_k}\rangle_{0,n,\beta}^X.$$
$\mathcal{F}^0_X$ is a formal power series in variables $t_{i,j}, 1\leq i\leq \text{rank}H^*(X, \mathbb{C}), j\geq 0$ with coefficients in $\mathbb{Q}[[\overline{NE}(X)]]$  and $Q^\beta$ is (again) a formal variable.
Using Theorem 2.3 we prove 
\begin{thm}\label{decomp_genus_0}
$$\mathcal{F}^0_{\clG}(\{t_{i\rho, j}\}_{1\leq i\leq \text{rank}H^*(X,\mathbb{C}), \rho\in \widehat{\mu_r}, j\geq 0}; Q)=\frac{1}{r^2}\sum_{[\rho]\in \widehat{\mu_r}}\mathcal{F}^0_X(\{t_{i\rho,j}\}_{1\leq i\leq \text{rank}H^*(X,\mathbb{C}), j\geq 0};Q_\rho),$$
where $Q_\rho$ is defined by the following rule: $$Q_\rho^\beta:=Q^\beta \chi_\rho\left(\exp\left(\frac{-2\pi i \int_\beta c_1(\mathcal{L})}{r}\right)\right),$$
and $\chi_\rho$ is the character associated to the representation $\rho$.
\end{thm}

Theorem \ref{decomp_genus_0} confirms the decomposition conjecture for genus $0$ Gromov-Witten theory of $\clG$.

\begin{rmk}
\begin{enumerate}
\item
If $X$ has generically semi-simple quantum cohomology, then an application of Givental's formula \cite{gi} shows that Theorem \ref{decomp_genus_0} implies the decomposition conjecture for $\clG$ in all genera.

\item
By definition a gerbe over $X$ is {\em essentially trivial} if it  becomes trivial after contracted product with the trivial $\mathcal{O}_X^*$-gerbe. For a finite abelian group $G$, any essentially trivial $G$-banded gerbe can be obtained as a fiber  product over the base of root gerbes. All the results presented here for root gerbes can be easily extended to this more general class of gerbes.

\end{enumerate}
\end{rmk}

\section{Results on toric gerbes}

Toric gerbes over toric orbifolds are toric Deligne-Mumford stacks in the sense of Borisov-Chen-Smith \cite{BCS}.  Any toric Deligne-Mumford stack can be constructed by taking a sequence of root of line bundles on the toric orbifolds, see \cite{JT}, \cite{FMN}.  

A toric Deligne-Mumford stack is defined in terms of a stacky fan $\mathbf{\Sigma}=(N,\Sigma,\beta)$, where $N$ is a finitely generated abelian group, $\Sigma\subset N_\mathbb{Q}=N\otimes_{\mathbb{Z}}\mathbb{Q}$ is a simplicial fan and $\beta: \mathbb{Z}^{n}\to N$ is a map determined by the elements $\{b_{1},\cdots,b_{n}\}$ in $N$.  By assumption, $\beta$ has finite cokernel and the images of $b_{i}$'s under the natural map $N\to N_\mathbb{Q}$ generate the simplicial fan $\Sigma$. The toric Deligne-Mumford stack $\XX(\mathbf{\Sigma})$  associated to $\mathbf{\Sigma}$ is defined to be  the quotient stack $[Z/G]$, where $Z$ is the open subvariety $\mathbb{C}^{n}\setminus\mathbb{V}(J_{\Sigma})$, $J_{\Sigma}$ is the irrelevant ideal of the fan, and $G$ is the product of an algebraic torus and a finite abelian group. The $G$-action on $Z$ is given by a group homomorphism $\alpha: G \to (\mathbb{C}^{*})^{n}$, where $\alpha$ is obtained by applying the functor  $\text{Hom}_{\mathbb{Z}}(-,\mathbb{C}^{*})$  to the Gale dual $\beta^{\vee}: \mathbb{Z}^{n}\to N^{\vee}$ of $\beta$ and $G=\text{Hom}_{\mathbb{Z}}(N^{\vee},\mathbb{C}^{*})$.

Every stacky fan $\mathbf{\Sigma}$ has an underlying {\em reduced} stacky fan $\mathbf{\Sigma_{red}}=(\overline{N},\Sigma,\overline{\beta})$, where $\overline{N}:=N/N_{tor}$, $\overline{\beta}: \mathbb{Z}^{n}\to \overline{N}$ is the natural projection given by the vectors $\{\overline{b}_{1},\cdots,\overline{b}_n\}\subseteq \overline{N}$. With these data one gets a toric Deligne-Mumford stack $\XX(\mathbf{\Sigma_{red}})=[Z/\overline{G}]$, where $\overline{G}=\text{Hom}_{\mathbb{Z}}(\overline{N}^{\vee},\mathbb{C}^{*})$ and $\overline{N}^{\vee}$ is the Gale dual $\overline{\beta}^{\vee}: \mathbb{Z}^{n}\to \overline{N}^{\vee}$ of the map $\overline{\beta}$. The stack $\XX(\mathbf{\Sigma_{red}})$ is a toric orbifold\footnote{I.e. the generic stabilizer is trivial.}, and can be obtained by rigidifying $\XX(\mathbf{\Sigma})$. We assume that $\XX(\mathbf{\Sigma})$ and $\XX(\mathbf{\Sigma_{red}})$ are semi-projective (see e.g. \cite{JT} for definition).

\subsection{Orbifold cohomology}
The Chen-Ruan orbifold cohomology ring of a toric Deligne-Mumford stack has been computed\footnote{Strictly speaking what's computed in \cite{BCS} is the orbifold Chow ring. The computation for Chen-Ruan orbifold cohomology ring is identical.} in \cite{BCS}. We recall the answer. Let $M=N^{*}$ be the dual of $N$. Let $\mathbb{C}[N]^{\mathbf{\Sigma}}$ be the group ring of $N$, i.e. $\mathbb{C}[N]^{\mathbf{\Sigma}}:=\bigoplus_{c\in N}\mathbb{C}y^{c}$, $y$ is the formal variable. Define the following multiplication 
\begin{equation}\label{product}
y^{c_{1}}\cdot y^{c_{2}}:=\begin{cases}y^{c_{1}+c_{2}}&\text{if
there is a cone}~ \sigma\in\Sigma ~\text{such that}~ \overline{c}_{1}, \overline{c}_{2}\in\sigma\,,\\
0&\text{otherwise}\,.\end{cases}
\end{equation}
Let $\mathcal{I}(\mathbf{\Sigma})$ be the ideal in
$\mathbb{C}[N]^{\mathbf{\Sigma}}$  generated by the elements $\sum_{i=1}^{n}\theta(b_{i})y^{b_{i}}, \theta\in M$. Then there is an isomorphism of $\mathbb{Q}$-graded algebras:
$$H_{CR}^{*}\left(\XX(\mathbf{\Sigma}), \mathbb{C} \right)\cong \frac{\mathbb{C}[N]^{\mathbf{\Sigma}}}{\mathcal{I}(\mathbf{\Sigma})}.$$

The natural map $N\to \overline{N}$ induces a map $\mathbf{\Sigma}\to\mathbf{\Sigma_{red}}$, which in turn induces a map of toric stacks $\XX(\mathbf{\Sigma})\to \XX(\mathbf{\Sigma_{red}})$. According to \cite{FMN}, \cite{JT}, this exhibits $\XX(\mathbf{\Sigma})$ as an $N_{tor}$-gerbe  over $\XX(\mathbf{\Sigma_{red}})$.  Moreover $\XX(\mathbf{\Sigma})$ is obtained from $\XX(\mathbf{\Sigma_{red}})$ as a tower of root gerbes.

Since $N=\overline{N}\oplus N_{tor}$, an element $c\in N$ has a unique decomposition $c=(\overline{c},\alpha)$ with $\overline{c} \in \overline{N}$ and $\alpha\in N_{tor}$. In particular we have $$b_i=(\overline{b}_i, \alpha_i)\in \overline{N}\oplus N_{tor}, 1\leq i\leq n.$$ This defines the elements $\alpha_i\in N_{tor}, 1\leq i\leq n$. 

Let $\widehat{N_{tor}}$ be the set of isomorphism classes of irreducible representations of 
$N_{tor}$. Since $N_{tor}$ is abelian, the set $\widehat{N_{tor}}$ is identified with the set of linear characters of $N_{tor}$. 

For $c\in N$, write $\overline{c}=\sum_{\overline{b}_{i}\subset \sigma(\overline{c})}a_{i}\overline{b}_{i}$, where 
$\sigma(\overline{c})$ is the minimal cone in $\Sigma$ containing 
$\overline{c}$. For $[\rho]\in \widehat{N_{tor}}$ denote by $\chi_\rho$ the associated linear character. Define
$$y^{\overline{c}, \rho}:=\frac{1}{|N_{tor}|}\left(
\sum_{\alpha\in N_{tor}}\chi_{\rho}(\alpha^{-1})\cdot y^{(\bar{c}, \alpha)}\right)\cdot \chi_{\rho}\left(\sum_{i=1}^n a_{i}\alpha_{i}\right)\in H_{CR}^*(\XX(\mathbf{\Sigma}), \mathbb{C}).$$

For each $[\rho]\in \widehat{N_{tor}}$, let $H_{CR}^*(\XX(\mathbf{\Sigma_{red}}))_{[\rho]}:=H_{CR}^*(\XX(\mathbf{\Sigma_{red}}), \mathbb{C})$. The direct sum $$\bigoplus_{[\rho]\in \widehat{N_{tor}}} H_{CR}^*(\XX(\mathbf{\Sigma_{red}}))_{[\rho]}$$ inherits a structure of a $\mathbb{Q}$-graded algebra from its summands. For $\bar{c}\in \overline{N}$ let $y_{\rho}^{\bar{c}}$ denote the element $y^{\bar{c}}$ in the summand $H_{CR}^*(\XX(\mathbf{\Sigma_{red}}))_{[\rho]}$ indexed by $[\rho]$.

\begin{thm}\label{orb_coh_isom}
The map 
$$\bigoplus_{[\rho]\in \widehat{N_{tor}}} H_{CR}^*(\XX(\mathbf{\Sigma_{red}}))_{[\rho]}\longrightarrow H_{CR}^*(\XX(\mathbf{\Sigma}), \mathbb{C}), \quad y_{\rho}^{\bar{c}}\mapsto y^{\bar{c}, \rho},$$
is an isomorphism of $\mathbb{Q}$-graded algebras.
\end{thm}
Theorem \ref{orb_coh_isom} is easily deduced as a corollary of the calculations in \cite{BCS}. 

\subsection{Gromov-Witten theory}
A detailed discussion of the basics of orbifold Gromov-Witten theory can be found in \cite{ts}. We obtain a comparison of Gromov-Witten theory of $\XX(\mathbf{\Sigma})$ and $\XX(\mathbf{\Sigma_{red}})$. Our result is most conveniently stated in terms of the {\em total descendant potential}, which is the generating function of all descendant Gromov-Witten invariants. The main tool is a detailed calculation of Gromov-Witten invariants of toric stacks \cite{CCIT2}.

\begin{thm}\label{decomp_toric}
The total descendant potential of $\XX(\mathbf{\Sigma})$ is a sum of $|N_{tor}|$ copies (indexed by $\widehat{N_{tor}}$) of the total descendant potential of $\XX(\mathbf{\Sigma_{red}})$, under the following change of variables:
\begin{enumerate}
\item
the cohomology variables are changed according to the isomorphism in Theorem \ref{orb_coh_isom};

\item 
the Novikov  variables in the descendant potential of $\XX(\mathbf{\Sigma_{red}})$ indexed by $[\rho]\in \widehat{N_{tor}}$ are rescaled as follows: $Q^d\mapsto Q^d\chi_\rho(\sum_{i=1}^na_i\alpha_i)$ for $d=\sum_{i=1}^na_ie_i$ in the Mori cone $\overline{\text{NE}}(\XX(\mathbf{\Sigma_{red}}))\subset Ker(\bar{\beta})\otimes_\mathbb{Z}\mathbb{R}\subset \mathbb{R}^n$; 

\item the genus variable $\hbar$ in the total descendant potential of $\XX(\mathbf{\Sigma_{red}})$ is rescaled by $1/|N_{tor}|$.

\end{enumerate}
\end{thm}

This Theorem verifies the decomposition conjecture for the Gromov-Witten theory of the gerbe $\XX(\mathbf{\Sigma})\to\XX(\mathbf{\Sigma_{red}})$. See \cite{AJT2} for more details.

\begin{rmk}
\begin{enumerate}
\item
Our results in fact are valid more generally for toric gerbes over toric Deligne-Mumford stacks which are not necessarily orbifolds. Details will be given in \cite{AJT2}.

\item
When the base is a $\mathbb{P}^1$-stack with at most two cyclic stack points, our results have also been proven by P. Johnson \cite{Johnson} by a completely different method.

\end{enumerate}
\end{rmk}

\begin{example}
We illustrate part of Theorem \ref{decomp_toric} concerning quantum cohomology rings in the example $\mathbb{P}(4,6)\to \mathbb{P}(2,3)$, which is the $\mu_2$-gerbe obtained as the stack of square roots of $\mathcal{O}_{\mathbb{P}(2,3)}(1)$. 

In \cite{AGV2} the quantum cohomology rings of $\mathbb{P}(4,6)$ and $\mathbb{P}(2,3)$ (and more generally all weighted projective lines) are computed:
\begin{equation*}
\begin{split}
QH_{orb}^*(\mathbb{P}(2,3), \mathbb{C})\simeq \mathbb{C}[[q]][x,y]/(xy-q, 2x^2-3y^3),\\
QH_{orb}^*(\mathbb{P}(4,6),\mathbb{C})\simeq \mathbb{C}[[q]][u,v,\xi]/(uv-q\xi, 2u^2\xi-3v^3, \xi^2-1).
\end{split}
\end{equation*}
\end{example}
For $i=0, 1$ let $QH_{orb}^*(\mathbb{P}(2,3), \mathbb{C})_i$ be a copy of $QH_{orb}^*(\mathbb{P}(2,3),\mathbb{C})$ with generators $x_i, y_i$ and $q$ rescaled by $(-1)^i$:
$$QH_{orb}^*(\mathbb{P}(2,3), \mathbb{C})_i=\mathbb{C}[[q]][x_i,y_i]/(x_iy_i-(-1)^iq, 2x_i^2-3y_i^3).$$
Let ${\bf 1}_0:=\frac{1}{2}(1+\xi), {\bf 1}_{1}:=\frac{1}{2}(1-\xi)$ and $u_i:=(-1)^iu{\bf 1}_i, v_i:=(-1)^iv{\bf 1}_i$. Then it is easy to check that the additive basis $\{{\bf 1}_i, u_i, v_i, v_i^2| i=0,1\}$ determines an isomorphism of algebras:
\begin{equation*}
\begin{split}
QH_{orb}^*(\mathbb{P}(4,6),\mathbb{C})\simeq QH_{orb}^*(\mathbb{P}(2,3), \mathbb{C})_0\oplus QH_{orb}^*(\mathbb{P}(2,3),\mathbb{C})_1,\\
{\bf 1}_i\mapsto 1\in QH_{orb}^*(\mathbb{P}(2,3), \mathbb{C})_i, u_i\mapsto x_i, v_i\mapsto y_i.
\end{split}
\end{equation*}

For instance, 
\begin{equation*}
\begin{split}
&{\bf 1}_0{\bf 1}_1=\frac{1}{4}(1-\xi^2)=0,\quad  {\bf 1}_0{\bf 1}_0=\frac{1}{2}(1+\xi)={\bf 1}_0,\quad  {\bf 1}_1{\bf 1}_1={\bf 1}_1,\\ 
&u_0v_1=0,\quad u_1v_0=0,\\
&u_0v_0=\frac{1}{2}(uv+uv\xi)=uv{\bf 1}_0=q\xi{\bf 1}_0=q{\bf 1}_0,\\
&u_1v_1= \frac{1}{2}(uv-uv\xi)=uv\frac{1}{2}(1-\xi)=q\xi\frac{1}{2}(1-\xi)=q\frac{1}{2}(\xi-1)=-q{\bf 1}_1,\\
&2u_i^2=2u^2{\bf 1}_i=3v^3\xi{\bf 1}_i=3v^3(-1)^i{\bf 1}_i=3v_i^3, \quad i=0,1.
\end{split}
\end{equation*}

\end{document}